\documentclass[11pt]{article}

\usepackage[active]{srcltx} % SRC Specials: DVI [Inverse] Search
\usepackage{amsmath,amsxtra,amssymb,amsthm}
\usepackage{latexsym,array,verbatim,fancybox,multirow}
\usepackage{epsfig,graphicx,url,setspace}
\usepackage{algorithm,algorithmic}

\usepackage{times,pslatex}
\usepackage{mathptmx}      % use Times fonts if available on your TeX system
\usepackage{ae,aecompl}

\newtheorem{prop}{Proposition}
\theoremstyle{definition}

\newtheorem{prob}{Problem}
\theoremstyle{remark}

\numberwithin{equation}{section}

\newcommand{\norm}[1]{\left\Vert#1\right\Vert}

\newcommand{\Real}{\mathbb R}

\newcommand{\mc}{\mathcal}

\begin{document}

\title{Dual Control with Active Learning \\ using Gaussian Process Regression} %\footnotemark[2]}

\author{Tansu Alpcan \\
              Technical University Berlin \\
	      Deutsche Telekom Laboratories \\
              \textit{alpcan@sec.t-labs.tu-berlin.de}         
}

\date{}
%\keywords{Dual control, system identification, Bayesian learning, regression, Gaussian processes, logistic map, inverted pendulum.}

\maketitle

% \footnotetext[2]{Please ensure that you use the most up to date
% class file,
% available from the ACS Home Page at\\
% \href{http://www3.interscience.wiley.com/journal/4508/home}{\texttt{http://www3.interscience.wiley.com/journal/4508/home}}}

%\vspace{-6pt}
%===============================================================================

\begin{abstract}
In many real world problems, control decisions have to be made with limited information. The controller may have no a priori (or even posteriori) data on the nonlinear system, except from a limited number of points that are obtained over time. This is either due to high cost of observation or the highly non-stationary nature of the system. The resulting conflict between information collection
(identification, exploration) and control (optimization, exploitation) necessitates an active learning approach for iteratively selecting the control actions which concurrently provide the data points for system identification. This paper presents a dual control approach where the information acquired at each control step is quantified using the entropy measure from information theory and serves as the training input to a state-of-the-art Gaussian process regression (Bayesian learning) method. The explicit quantification of the information obtained from each data point allows for iterative optimization of both identification and control objectives. The approach developed is illustrated with two examples: control of logistic map as a chaotic system and position control of a cart with inverted pendulum.
\end{abstract}

% -----------------------------------------------------------
\section{Introduction} \label{sec:intro}

In many real world problems, control decisions have to be made with limited information. Obtaining extensive and 
accurate information about the controlled system can often be a costly and time consuming process. In some cases, acquiring detailed information on system characteristics may be simply infeasible due to high observation costs. In others, the observed system may be so nonstationary that by the time the information is obtained, it is already outdated due to system's fast-changing nature. Therefore, the only option left to the controller is to develop a strategy for collecting information efficiently and choose a model to estimate the ``missing portions'' of the system in order to control it according to a given objective.

A variant of this problem has been well-known in the control literature since 1960s as \textit{dual control}. The underlying concept in dual control is obtaining good process information through perturbation while controlling it. The controller has necessarily dual goals. First the controller must control the process as well as possible. Second, the controller must inject a probing  signal or perturbation to get more information about the process. By gaining more process information better control can be achieved in the future~\cite{wittenmark1}.

The problem considered here differs from the classical dual control problem in the very limited amount of information available to the controller. The controller here cannot aim to identify the system first to obtain better performance in the future due to non-stationarity and/or prohibitive observation costs. Furthermore, the perturbation idea is not fully applicable since each action-observation pair provides a single data point for identifying the nonlinear discrete-time system, unlike in the identification of (linear) continuous-time systems. %, perturbation does not play a significant role here.

This paper approaches the ``dual control'' problem from a Bayesian perspective. Gaussian processes (GP) are utilized as a state-of-the-art regression (function estimation) method for identifying the underlying state-space equations of the discrete-time nonlinear
system from observed (training) data. More importantly, the adopted GP (Bayesian) framework allows explicit quantification of information, which each observed data point provides within the a-priori chosen model. Hence, the information collection goal
can be explicitly combined with the control objectives and posed as a (weighted-sum, multi-objective) optimization problem
based on one (or multi-) step lookahead. This results in a joint and iterative scheme of active learning and control.

\begin{figure}[htp]
 \centering
 \includegraphics[width=0.9\columnwidth]{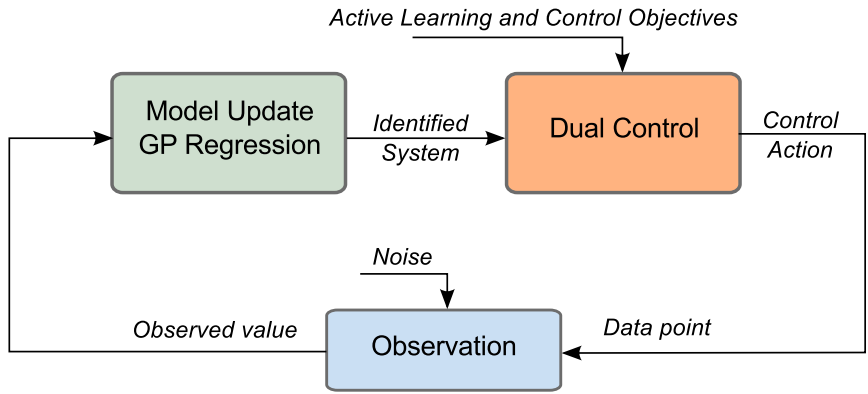}
 \caption{The underlying model of the dual control approach.}
 \label{fig:model1}
\end{figure}

The proposed approach consists of three main parts: observation, update of GP for regression, and optimization to determine the next control action. These three steps, shown in Figure~\ref{fig:model1} are taken iteratively to achieve the dual objectives of identification and control. 
%As a result of its iterative nature, this approach can be considered in a sense similar to the well-known Expectation-Maximization algorithm \cite{Bishopbook}.

Observations, given that they are a scarce resource in the class of problems considered, play an important role in this approach. Uncertainties in the observed quantities can be modeled as additive noise. Likewise, properties (variance or bias) of additive noise can be used to model the reliability of (and bias in) the data points observed. GPs provide a straightforward mathematical structure for incorporating these aspects to the model under some simplifying assumptions.

The set of observations collected provide the (supervised) training data for GP regression in order to estimate the characteristics of the function or system at hand. This process relies on the GP methods, which will be described in Subsection~\ref{sec:gp}. Thus, at each iteration an up-to-date description of the function or system is obtained
based on the latest observations. 

% Specifically, $\hat f$ provides an estimate of the original function $f$.\footnote{See \cite[Chap 7.2]{GPbook} for a discussion on asymptotic analysis of GP regression. It should not be noted, however, that asymptotic properties are of little relevance to the problem at hand. }
% Assuming an additive Gaussian noise model, the noise variance $\sigma$ can be used to model uncertainties, e.g. older and noisy data resulting in higher $\sigma$ values.

The final step of the approach provides a basis for determining the next control action based on an optimization process that takes into account dual objectives. The information measurement aspect of these objectives will be discussed in Subsection~\ref{sec:obsinfo}. An important issue here is the fact that there are infinitely many candidate points in this optimization process, but in practice only a finite collection of them can be evaluated. 

The investigated approach incorporates many concepts that have been implicitly considered by heuristic schemes, and  builds upon  results from seemingly disjoint but relevant fields such as information theory, machine learning, optimization, and control theory.  Specifically, it combines concepts from these fields by
% The proposed framework brings together the mentioned fields and combines these concepts by
\begin{itemize}
 \item explicitly quantifying the information acquired using the entropy measure from information theory,
 \item modeling and estimating the (nonlinear) controlled system adopting a Bayesian approach and using Gaussian processes as a state-of-the-art regression method, 
 \item using an iterative scheme %similar to Expectation-Maximization (EM) 
for observation, learning, and control,
 \item capturing all of these aspects under the umbrella of a multi-objective ``meta'' optimization and control formulation.
\end{itemize}

Despite methods and approaches from machine (statistical) learning are heavily utilized in this framework, the problem at hand is very different from many classical machine learning ones, even in its learning aspect. In most classical application domains of
machine learning such as data mining, computer vision, or image and voice recognition, the difficulty is often in handling significant amount of data in contrast to lack of it. Many methods such as Expectation-Maximization (EM) inherently make this assumption, except from ``active learning'' schemes \cite{Bishopbook}. Information
theory plays plays an important role in evaluating scarce (and expensive) data and developing strategies for obtaining it. Interestingly, data scarcity converts at the same time the disadvantages of some methods into advantages, e.g. the scalability problem of Gaussian processes.

It is worth noting that the class of problems described here are much more frequently encountered in practice than it may first seem. Social systems and economics, where information is scarce and systems are very non-stationary by nature constitute an
important application domain. The control framework proposed is further applicable to a wide variety of fields due to its fundamentally adaptive nature. One example is decentralized resource allocation decisions in networked and complex systems, e.g. wired and wireless networks, where parameters change quickly and global information on network characteristics is not available at the local decision-making nodes. Another example is security and information technology risk management in large-scale organizations, where acquiring information on individual subsystems and processes can be very costly. Yet another example application is in biological systems where individual organisms or subsystems operate autonomously (even if they are part of a larger system) under limited local information.

% 
% -----------------------------------------------------------
\section{Methodology} \label{sec:methods}

This section summarizes the results in \cite{valuetools11} and presents the underlying methods that are utilized within the dual control framework. First, the regression model and Gaussian Processes (GP) are presented. Subsequently, modeling and measurement of information is discussed using (Shannon) information theory.
%Finally, an integrated model combining these concepts and methods is introduced.

\subsection{Regression and Gaussian Processes (GP)} \label{sec:gp}
The system identification problem here involves inferring the nonlinear function(s) $f$ in the state-space equations 
describing the system using the set of observed data points. This is known as \textit{regression} in machine learning literature, which is a \textit{supervised learning} method since the data observed here is at the same time the training data. This learning process involves selection of a ``model'', where the learned function $\hat f$ is, for example, expressed in terms of a set of parameters and specific basis functions. Gaussian processes (GP) provide a nonparametric alternative to this but follow in spirit the same idea. 

The main goal of regression involves a trade-off. On the one hand, it tries to minimize the \textit{observed} error between $f$ and $\hat f$. On the other, it tries to infer the ``real'' shape of $f$ and make good estimates using 
$\hat f$ even at unobserved points (generalization). If the former is overly emphasized, then one ends up with ``over fitting'', which means $\hat f$ follows $f$ closely at observed points but has weak predictive value at unobserved ones. This delicate balance is usually achieved by balancing the prior ``beliefs'' on the nature of the function, captured by the model (basis functions), and fitting the model to the observed data. 

This paper focuses on Gaussian Process \cite{GPbook} as the chosen regression method within the proposed dual control approach without loss of any generality. There are multiple reasons behind this preference. Firstly, GP provides an elegant mathematical method for easily combining many aspects of the approach. Secondly, being a nonparametric method GP eliminates any discussion on model degree. Thirdly, it is easy to implement and understand as it is based on well-known Gaussian probability concepts. Fourthly, noise in observations is immediately taken into account if it is modeled as Gaussian. Finally, one of the main drawbacks of GP namely being computational heavy, does not really apply to the problem at hand since the amount of data available is already very limited.

It is not possible to present here a comprehensive treatment of GP. Therefore, a very rudimentary overview is provided next within the context of the control problem. Consider a set of $M$ data points 
$$\mc D=\{x_1, \ldots, x_M\},$$
where each $x_i \in \mc X$ is a $d-$dimensional vector, and the corresponding vector of scalar values is $f(x_i), \; i=1,\ldots,M$. Assume that the observations are distorted by a zero-mean Gaussian noise, $n$ with variance $\sigma \sim \mc N(0,\sigma)$. Then, the resulting observations is a vector of Gaussian $y=f(x)+n \sim \mc N(f(x),\sigma)$. 

A GP is formally defined as a collection of random variables, any finite number of which have a joint Gaussian distribution \cite{GPbook}. It is completely specified by its mean function $m(x)$ and covariance function $C(x,\tilde x)$, where
$$ m(x)=E[\hat f(x)] $$
and 
$$C(x,\tilde x)=E[(\hat f(x)-m(x))(\hat f(\tilde x)-m(\tilde x))], 
\; \forall x, \tilde x \in \mc D. $$

Let us for simplicity choose $m(x)=0$. Then, the GP is characterized entirely by its covariance function $C(x,\tilde x)$. Since the noise in observation vector $y$ is also Gaussian, the covariance function can be defined as the sum of a \textit{kernel function} $Q (x,\tilde x)$ and the diagonal noise variance 
\begin{equation} \label{e:gcov}
 C(x,\tilde x) = Q (x,\tilde x) + \sigma I, \; \forall \, x, \tilde x \in \mc D ,
\end{equation}
where $I$ is the identity matrix. While it is possible to choose here any (positive definite) kernel $Q(\cdot,\cdot)$, one classical choice is 
\begin{equation} \label{e:gaussiankernel}
 Q(x,\tilde x)=\exp \left[-\frac{1}{2}\norm{x -\tilde x}^2 \right].
\end{equation}
Note that GP makes use of the well-known \textit{kernel trick} here by representing an infinite dimensional continuous function
using a (finite) set of continuous basis functions and associated vector of real parameters in 
accordance with the \textit{representer theorem} \cite{schoelkopfbook}.

The (noisy)\footnote{The special case of perfect observation without noise is handled the same way as long as the kernel function $Q(\cdot,\cdot)$ is positive definite.} training set $(\mc D, y)$ is used to define the corresponding GP, $\mc{GP} (0,C(\mc D))$, through the $M \times M$ covariance function $C(\mc D)=Q+\sigma I$, where the conditional Gaussian distribution of any point outside the training set, $\bar y \in \mc X, \bar y \notin \mc D$, given the training data $(\mc D, t)$  can be computed as follows. Define the vector 
\begin{equation} \label{e:k}
 k(\bar x)=[Q(x_1,\bar x), \ldots Q(x_M,\bar x)]
\end{equation}
and scalar 
\begin{equation} \label{e:kappa}
\kappa=Q(\bar x,\bar x)+\sigma.
\end{equation}
Then, the conditional distribution  $p(\bar y | y)$ that characterizes the $\mc{GP} (0,C)$ is a Gaussian $\mc N(\hat f,v)$ with mean $\hat f$ and variance $v$,
\begin{equation} \label{e:gp1}
 \hat f(\bar x)=k^T C^{-1} y \text{ and } v(\bar x)=\kappa - k^T C^{-1} k .
\end{equation}

This is a key result that defines GP regression as the mean function $\hat f(x)$ of the Gaussian distribution and provides a prediction of the function $f(x)$. At the same time, it belongs to the well-defined class $\hat f \in\mc F$, which is the set of all possible sample functions of the GP
$$\mc F := \{\hat f(x): \mc X \rightarrow \Real \text{ such that }  \hat f \in \mc{GP} (0,C(\mc D)),\; \forall \mc D, \, C \} ,$$
where $ C(\mc D)$ is defined in (\ref{e:gcov}) and $\mc{GP}$ through (\ref{e:k}), (\ref{e:kappa}), and (\ref{e:gp1}), above.
Furthermore, the variance function $v(x)$ can be used to measure the uncertainty level of the predictions provided by $\hat f$, which will be discussed in the next subsection.

%\textbf{TODO}: A simple numerical example to illustrate GP

\subsection{Quantifying the Information in Observations} \label{sec:obsinfo}

Each observation provides a data point to the regression problem (estimating $f$ by constructing $\hat f$) as discussed in the previous subsection. %Many works in the learning literature consider the ``training'' data used in regression available (all at once or sequentially) and do not discuss the possibility of the decision maker influencing or even optimizing the data collection process. 
\textit{Active learning}  addresses the question of ``how to quantify information obtained and optimize the observation process?''. Following
the approach discussed in \cite{MacKaydataselect,MacKaybook}, the approach here provides a precise answer to this question.

Making any decision on the next (set of) observations in a principled manner necessitates first \textit{measuring the information obtained from each observation within the adopted model}. It is important to note that the information measure here is dependent on the chosen model. For example, the same observation provides a different amount of information to a random search model than a GP one.

Shannon information theory readily provides the necessary mathematical framework for measuring the information content of a variable. Let $p$ be a probability distribution over the set of possible values of a discrete random variable $A$. The \textbf{entropy} of the random variable is given by
$H(A)=\sum_i p_i \log_2 (1/p_i)$, which quantifies the amount of uncertainty. Then, the information obtained from an observation on the variable, i.e. reduction in uncertainty, can be quantified simply by taking the difference of its initial and final entropy, 
$$\mc I=H_0 - H_1. $$ 

It is important here to avoid the common conceptual pitfall of equating entropy to information itself as it is sometimes done in communication theory literature. Since this issue is not of great importance for the class of problems considered in communication theory, it is often ignored. However, the difference is of conceptual importance in this 
problem.\footnote{See \url{http://www.ccrnp.ncifcrf.gov/~toms/information.is.not.uncertainty.html} for a detailed discussion.}
In this case, (Shannon) \textit{information is defined as a measure of the decrease of uncertainty after (each) observation (within a given model)}. %This can be best explained with the following simple example \cite{valuetools11}.

%\subsection{Quantifying Information in GP}

To apply this idea to GP, let the zero-mean multivariate Gaussian (normal) probability distribution be denoted as
\begin{equation} \label{e:multivargauss}
 p(x)=\dfrac{1}{\sqrt{2\pi |C_p(x)}|} \exp \left( -\frac{1}{2}[x-m]^T|C_p(x)|^{-1} [x-m]\right),
\end{equation}
where $x \in \mc X$, $|\cdot|$ is the determinant, $m$ is the mean (vector) as defined in (\ref{e:gp1}), and $C_p(x)$ is the covariance matrix as a function of the newly observed point $x \in \mc X$ given by
\begin{equation} \label{e:covx}
 C_p(x)=\left[ 
\begin{array}{cccc}
 &   &  &  \\
 & C(\mc D) &  & k(x)^T \\
 &   &  &  \\
 & k(x) &  & \kappa
\end{array}
\right] .
\end{equation}
Here, the vector $k(x)$ is defined in (\ref{e:k}) and $\kappa$ in (\ref{e:kappa}), respectively. The matrix $C(\mc D)$ is the covariance matrix based on the training data $\mc D$ as defined in (\ref{e:gcov}). 

The entropy of the multivariate Gaussian distribution (\ref{e:multivargauss}) is \cite{entropygaussian} 
$$ H(x)=\dfrac{d}{2}+\dfrac{d}{2}\ln(2\pi)+\dfrac{1}{2} \ln |C_p(x)| ,$$
where $d$ is the dimension. Note that, this is the entropy of the GP estimate at the point $x$ based on the available data $\mc D$. The aggregate entropy of the function on the region $\mc X$ is given by
\begin{equation} \label{e:aggentropy}
 H^{agg}:=\int_{x \in \mc X} \dfrac{1}{2} \ln |C_p(x)| dx.
\end{equation}

The problem of choosing  a new data point $\hat x$ such that the information obtained from it within the
GP regression model is maximized can be formulated as:
\begin{eqnarray} \label{e:infocollect1}
 \hat x=\arg \max_{\tilde x} \mc I= \arg \max_{\tilde x} \int_{x \in \mc X} \left[ H_0 - H_1 \right] \, dx \\
\nonumber = \arg \min_{\tilde x}  \int_{x \in \mc X} \dfrac{1}{2} \ln |C_q(x,\tilde x)| dx,
\end{eqnarray}
where the integral is computed over all $x \in \mc X$, and the covariance matrix $C_q(x, \tilde x)$ is defined as
\begin{equation} \label{e:covxbar}
 C_q(x, \tilde x)=\left[ 
\begin{array}{ccccc}
 &   &         &          &     \\
 & C(\mc D)&   & k^T(\tilde x) & k^T(x) \\
 &   &         &          &  \\
 &  k(\tilde x) &   & \tilde \kappa & Q(x,\tilde x) \\
 &  k(x) &  & Q(x,\tilde x) & \kappa
\end{array}
\right] ,
\end{equation}
and $\tilde \kappa=Q(\tilde x,\tilde x)+\sigma$. Here, $C(\mc D)$ is a $M \times M$ matrix and $C_q$ is a $(M+2) \times (M+2)$ one, whereas $\kappa$ and  $Q(x,\tilde x)$ are scalars and $k$ is a $M \times 1$ vector. 
% Thus, the problem of choosing the next ``most informative'' data point $\tilde x$ is equivalent to:
% $$ \max_{\tilde x} \mc I=\max_{\tilde x} \int_{x \in \mc X} \ln (1/|C_q(x,\tilde x)|) 
% %\equiv  \min_{\tilde x} \prod_x |C_q(x, \tilde x)| 
% .$$
This result from \cite{valuetools11} is summarized in the following proposition.

\begin{prop} \label{thm:GPsearch}
As a maximum information data collection strategy for a Gaussian Process with a covariance matrix $C(\mc D)$, the next observation $\hat x$ should be chosen in such a way that 
$$ \hat x= \arg \max_{\tilde x} \mc I= \arg  \min_{\tilde x}  \int_{x \in \mc X} \ln |C_q(x,\tilde x)| dx,$$
%$$ \tilde x= \arg \min_{\tilde x} \prod_x |C_q(x, \tilde x)|, $$
where $C_q(x, \tilde x)$ is defined in (\ref{e:covxbar}).
\end{prop}

\subsubsection*{An Approximate Solution to Information Maximization}

When making a decision on the next action through multi-objective optimization, there are (infinitely) many candidate points. A pragmatic solution to the problem of finding solution candidates is to (adaptively) sample the problem domain $\mc X$ to obtain the set  
$$\Theta:=\{x_1, \ldots, x_T : x_i \in \mc X, \, x_i \notin \mc D, \; \forall i \}$$ 
that does not overlap with known points. In low (one or two) dimensions, this can be easily achieved through grid
sampling methods. In higher dimensions, (Quasi) Monte Carlo schemes can be utilized. For large problem domains, the current domain of interest $\mc X$ can be defined around the last or most promising observation in such a way that such a sampling is computationally feasible. 
Likewise, multi-resolution schemes can also be deployed to increase computational efficiency.

Given a set of (candidate) points $\Theta$ sampled from $\mc X$, the result in Proposition~\ref{thm:GPsearch} can be revisited. The problem in (\ref{e:infocollect1}) is then approximated  \cite{tempobook} by
\begin{eqnarray} \label{e:infocollect2}
 \max_{\tilde x} \mc I \approx \min_{\tilde x} \sum_{x \in \Theta} \ln |C_q(x,\tilde x)| \\
 \nonumber \Rightarrow  \hat x= \arg \min_{\tilde x \in \Theta}  \prod_{x \in \Theta} |C_q(x, \tilde x)|,
\end{eqnarray}
using monotonicity property of the natural logarithm and the fact that the determinant of a covariance matrix is non-negative.
%where $C_q(x, \tilde x)$ is defined in (\ref{e:covxbar}). 
Thus, the following counterpart of Proposition~\ref{thm:GPsearch}
is obtained:
\begin{prop} \label{thm:GPsearch2}
As an approximately maximum information data collection strategy for a Gaussian Process with a covariance matrix $C(\mc D)$ and given a collection of candidate points $\Theta$, the next observation $\hat x \in \Theta$ should be chosen in such a way that 
$$ \hat x=  \arg  \min_{\tilde x \in \Theta}  \prod_{x \in \Theta} |C_q(x, \tilde x)| \approx \arg \max_{\tilde x \in \Theta} \mc I,$$
%$$ \tilde x= \arg \min_{\tilde x} \prod_x |C_q(x, \tilde x)|, $$
where $C_q(x, \tilde x)$ is given in (\ref{e:covxbar}).
\end{prop}

Although it is an approximation, finding a solution to the optimization problem in Proposition~\ref{thm:GPsearch2} can still be computationally  costly. Therefore, a greedy algorithm is proposed as a computationally simpler alternative.
% Let,  $x^* \in \Theta$ be defined as 
% $$ x^* :=\arg \max_{x \in \Theta} |C_p(x)|=|C(\mc D)|\, |\kappa(x) - k(x) C^{-1}(\mc D) k^T(x) |,$$
% where the matrix $C_p$ is given by (\ref{e:covx}) \cite{matrixcookbook}.  The first term above, $|C(\mc D)|$ is fixed and the second one, 
% $$|\kappa(x) - k(x) C^{-1}(\mc D) k^T(x) |, $$
% is the same as the GP variance $v(x)$ in (\ref{e:gp1}). Hence, the sample $x^*$ is one of those with the maximum variance in the set $\Theta$, given current data $\mc D$. 
% 
% It follows from (\ref{e:covxbar}) and basic matrix theory that if $\tilde x=x$ for a given $x$ then
% $ |C_q(x, \tilde x)|$ is minimized. As a simplification, ignore the dependencies between $C_q(x, \tilde x)$ matrices for different $x \in \Theta$. Then, 
Choosing the maximum variance $\hat x$ as
$$ \hat x = \arg \max_{\tilde x \in \Theta} v(\tilde x) \approx \arg \min_{\tilde x \in \Theta}  \prod_{x \in \Theta} |C_q(x, \tilde x)|,$$
leads to a large (possibly largest) reduction in $\prod_{x \in \Theta} |C_q(x, \hat x)|$, and hence
provides a rough approximate solution to (\ref{e:infocollect2}) and to the result in Proposition~\ref{thm:GPsearch}.
This result from \cite{valuetools11} is consistent with widely-known heuristics such as ``maximum entropy'' or ``minimum variance'' methods \cite{activelearning} and a variant has been discussed in \cite{MacKaydataselect}.

\begin{prop} \label{thm:GPsearch3}
Given a Gaussian Process with a covariance matrix $C(\mc D)$ and a collection of candidate points $\Theta$, an approximate solution to the maximum information data collection problem defined in Proposition~\ref{thm:GPsearch} is to choose
the sample point(s) $\tilde x$ in such a way that it has (they have) the maximum variance within the set $\Theta$.
\end{prop}

% -----------------------------------------------------------
\section{Dual Control with Limited Information} \label{sec:dynamic}

% The previous sections have focused mainly on static optimization. However, the problem of controlling dynamic systems with limited information is at least as important and widely applicable as its static counterpart. Fortunately, the methodology presented in Section~\ref{sec:methods} is directly applicable to dynamic control with limited information by replacing \textit{regression} with \textit{system identification}. System identification involves mathematical description and modeling of dynamical systems from observation and measured data. Although the same underlying principles of the framework are applicable in this case, the dynamic nature of the problem also leads to minor changes. One of the main differences is the fact that the system states are now influenced indirectly through control actions. The data points used for identifying the underlying system mapping can only be selected indirectly (unlike static optimization) and under the constraints imposed by the nature of the ``control'' in the dynamic system at hand.

Consider a nonlinear discrete-time representation of a dynamical system that evolves on a $d-$dimensional state space $\mc X^d \subset \Real^{d}$ steered by control actions chosen from an $e-$dimensional space $\mc U^e \subset \Real^e$. Usually, the dimension of the control space is smaller than the state one, $e \leq d$. It is assumed here for simplicity that both control and state spaces are nonempty, convex, and compact.
%It is also possible to choose infinite versions of both spaces in certain cases without affecting the main aspects of the analysis. 
The system states evolve according to
\begin{equation} \label{e:dyn1}
  x_i(t+1)= f_i (x(t),u(t)), \;\; i=1,\ldots,d \, ,
\end{equation}
where $x(t) \in \mc X^{d}$, $x_i(t)$ is a scalar, $t=1,\ldots$ denotes discrete time instances, and each 
$f_i: \mc X^d \times \mc U^e  \rightarrow \Real$ 
is a possibly nonlinear function. States of dynamical systems are, however, often not observable. Therefore, define a mapping from the states to observable quantities $y$ as
\begin{equation} \label{e:dyn2}
  y_j(t)= g_j (x(t)), \;\; j=1,\ldots, \bar d \, ,
\end{equation}
where each $g_j: \mc X^d  \rightarrow \Real$ is possibly a nonlinear function, and $\bar d \leq d$.

If nothing is known about the dynamic system defined by (\ref{e:dyn1})-(\ref{e:dyn1}) in the beginning, and there is no observation or system noise, then the system can be simplified to its input-output relationship:
\begin{equation} \label{e:dyn3}
\begin{array}{ll}
 & y_j(t+1)=g_j\left(f[g^{-1}(y(t)),u(t)]\right) \\ \\
 \Rightarrow & y_j(t+1)=h_j \left(y(t),u(t)\right),\; \quad j=1,\ldots,\bar d \, , 
\end{array}
\end{equation}
where each $h_j: \mc X^{\bar d} \times \mc U^e \rightarrow \Real$ is possibly a nonlinear function. As a simplification, system and observation noise can be modeled as zero-mean Gaussian\footnote{Biased Gaussian noise can be easily handled by GPs by introducing a mean function, which we omit in this paper for simplicity.}. Thus, a noisy variant of system (\ref{e:dyn3}) is
\begin{equation} \label{e:dyn4}
  y_j(t+1)= h_j (y(t),u(t)) + n(t), \;\quad j=1,\ldots, \bar d \, ,
\end{equation}
where $n(t) \sim \mc N(0,\sigma)$ and $\sigma$ is the respective noise variance.

\subsection{Problem Formulation}

The dual control problem is defined as follows. Consider an unknown nonlinear discrete-time dynamic system, which has a control input and a (partially) observable output that is possibly distorted by noise. The control input may affect the
system linearly, which leads to a simpler problem, or its effect may be nonlinear and unknown to the decision maker. The objective of the decision maker is to control the system in such a way that it follows a given reference signal. Each action taken is
assumed to be very costly and the decision maker may only have limited time to satisfy dual goals of identification
and control. \textit{What is the best strategy to address this problem}?

Based on the discussion above, the described problem can be formulated more concretely. Let $r(t) \in \mc X^{\bar d}\; \forall t$ denote the $\bar d-$dimensional reference signal. The discrete-time nonlinear system can be modeled using (\ref{e:dyn4}), where $y(t)$ is the output, $u(t)$ is the control action, and $n(t)$ is the observation noise at time $t$.
Then, the following dual control problem is formulated.

\begin{prob}\label{prob:control1} 
[\textit{Dual Control under Limited Information}]
Let a discrete-time system be described by the following input-output relationship
$$  y_j(t+1)= h_j (y(t),u(t)) + n(t), \; \quad j=1,\ldots,\bar d$$
where $y(t)$ is the $\bar d-$dimensional output, $u(t)$ is the $e-$dimensional control action, and $n(t)\sim \mc N(0,\sigma)$ is a zero-mean Gaussian observation noise with variance $\sigma$ at time $t$. The function $h_j: \mc X^{\bar d} \times \mc U^e \rightarrow \Real$ is possibly nonlinear for all $j$. Given a $\bar d-$dimensional desired reference signal $r(t)$, what is the best control strategy (series of control actions) $\mu(t)$ such that 
$$ \mu(t)= \arg \min_{u(t)} \norm{y(t)-r(t)}, \; \forall t=1,\ldots, $$
$\norm{\cdot,\cdot}$ is a norm quantifying the mismatch between the observed and desired outputs?
\end{prob}

% The simplest (both conceptually and computationally) strategy to solve Problem~\ref{prob:control1} is random search on the control and state space. As such no attempt is made to ``learn'' the properties of the system (\ref{e:dyn4}). Unless, the system at hand is ``algorithmically random'' [\textbf{CITE}], which is rarely the case, this strategy is clearly suboptimal as it wastes the information collected on the system. 

If there was more information on the system available or more time for experimentation, one could have resorted to the rich literature on adaptive and robust control to find a solution. However, Problem~\ref{prob:control1}  differentiates from the ones in the classical adaptive and robust control literature by the fact that the decision maker starts with zero or very little prior information and a solution has to be found online while learning the system. This puts special emphasis on observations and quantifying information using the methods described in Section~\ref{sec:obsinfo}. 

Using GP regression for estimating the system dynamics
in (\ref{e:dyn4}) and Shannon information theory to measure and maximize the amount of information obtained with each observation, a model-based variant of Problem~\ref{prob:control1} is defined.

\begin{prob} \label{prob:control2}
[\textit{Model-based Control under Limited Information}]
Let a discrete time dynamic system be described by the following input-output relationship
$$  y_j(t+1)= h_j (y(t),u(t)) + n(t), \; \quad j=1,\ldots,\bar d$$
where $y(t)$ is the $\bar d-$dimensional output, $u(t)$ is the $e-$dimensional control action, and $n(t)\sim \mc N(0,\sigma)$ is a zero-mean Gaussian observation noise with variance $\sigma$ at time $t$. The function $h_j: \mc X^{\bar d} \times \mc U^e  \rightarrow \Real$ is possibly nonlinear for all $j$. The goal is to control the system in such a way that the output $y(t)$ follows a given $\bar d-$dimensional reference signal $r(t)$.

Let $\hat h$ be an estimate of system dynamics $h$ based on an a priori model and a set of observations. What is the best control strategy (series of control actions) $\mu(t)$ that solves the multi-objective problem with the following components?
\begin{itemize}
 \item \textit{Objective 1:} $ \min_{u \in \mc U}  \norm{y(t)-r(t)}$
 \item \textit{Objective 2:} $ \max_{u \in \mc U}  \mc I(\hat h,u(t))$
\end{itemize}
\end{prob}

The main (first) objective of Problem~\ref{prob:control2} is naturally the same as the one of Problem~\ref{prob:control1}. The second objective states the ``exploration'' or information collection aspect. 

As a side note, unlike the static optimization problem in \cite{valuetools11}, how close the estimated system dynamic $\hat h$ approximates the original one is not set as an objective. The reason behind this is the fact that the data points used for identifying $\hat h$ can only be selected indirectly through control actions $u$. Therefore, a reasonably complete identification of the system dynamics may be too costly. A partial identification relevant to the main objective is sufficient for the purpose here.

\subsection{Solution Approach}

\begin{figure}[htp]
 \centering
 \includegraphics[width=\columnwidth]{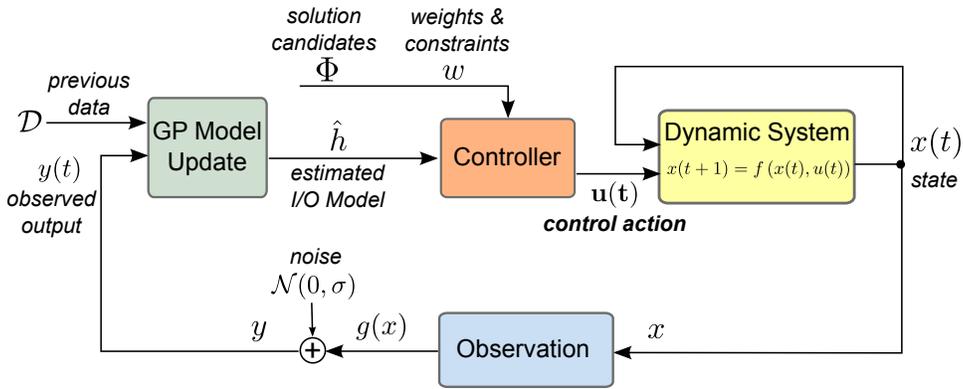}
 \caption{The dual control framework for identifying and controlling an unknown dynamic system with limited information.}
 \label{fig:model3}
\end{figure}

The solution approach to Problem~\ref{prob:control2} utilizes the methodology in Section~\ref{sec:methods}. The GP variance maximization approximates here the information maximization objective. A (random or grid-based) sampling scheme is adopted again for evaluating candidate solutions, in this case, a combination of the observed current state and available control actions. A weighted-sum scheme is utilized to combine the two objectives in Problem~\ref{prob:control2}. A visual depiction of the control framework is shown in Figure~\ref{fig:model3}. 

Since the problem is by its very nature iterative, the best control strategy has to be evaluated at the current state, taking into account newly received information and using the latest update of estimated system dynamics. As a starting point, a gradient or greedy algorithm is proposed which aims to balance both exploration and exploitation objectives. 

% There are, at the same time, practical and important differences between the dynamic and static versions of the decision making problem. The dynamic system is by its very nature iterative. Therefore, the best control strategy has to be evaluated at the current state, taking into account newly received information and using the latest update of estimated system dynamics. As a starting point, a gradient or greedy algorithm is proposed which aims to balance both exploration and exploitation objectives. 
 
\begin{prop} \label{prob:control3}
Let a discrete time dynamic system be described by the following input-output relationship
$$  y_j(t+1)= h_j (y(t),u(t)) + n(t), \; y(t) \in  \mc X^{\bar d},\; u(t) \in \mc U^e,   $$
$j=1,\ldots,\bar d$, where $n(t)\sim \mc N(0,\sigma)$ is a zero-mean Gaussian observation noise with variance $\sigma$. Further let $\Phi$ be a grid-based or randomly sampled set of available control actions $u$ from the control space $\mc U$. Given a reference signal $r(t) \in \mc X^{\bar d}$, define the optimization problem
\begin{equation} \label{e:contweigth1}
\begin{array}{c}
  \min_{u(t) \in \Phi} J(u,y,r,w)  \\ \\
 \Rightarrow \min_{u(t) \in \Phi} w_1 \, \norm{\hat y(t+1)-r(t)} - w_2 \, v(\hat y(t+1),\mu(t)), 
\end{array}
\end{equation}
where 
$$\hat y_j(t+1)=\hat h_j (y(t),u(t)) + n(t)$$ 
is the next estimated output using a GP based on control $u(t)$, and $v(\hat y(t+1),u(t))$ is the variance of the associated Gaussian as defined in (\ref{e:gp1}). 
The solution to this problem 
$$\mu(t)=\arg \min_{u(t) \in \Phi} J(u,y,r,w), \; t=1,\ldots $$ 
approximates the best control strategy under limited information, and hence approximately solves Problem~\ref{prob:control2}.
\end{prop}

Couple of remarks should be made at this point regarding the solution approach presented. Firstly, the approach in Proposition~\ref{prob:control3} constitutes a greedy one, which aims to solve the problem in shortest time based on available information and goes in the direction of the steepest gradient (here of the weighted sum of objectives).
The main concern here is whether such an algorithm gets stuck in a local minimum. This issue can be remedied at 
least partially by putting a higher weight to the information collection objective. 
Secondly, it is implicitly assumed here that the system at hand is at least partially observable and controllable. It is naturally difficult, if not impossible, to check such properties of an unknown system. Thus, the approach here can be interpreted also as a ``best effort'' one, which aims to achieve the best performance possible given controllability and observability limitations.

A summary of the solution approach discussed above for a specific set of choices is provided by Algorithm~\ref{alg:algctrl1}.
\begin{algorithm}[htbp]
   \caption{Dual Control with Limited Information}
   \label{alg:algctrl1}
\begin{algorithmic}
  \STATE {\bfseries Input:} Problem domain, GP meta-parameters, objective weights $[w_1, w_2]$, initial data set $\mc D$, reference signal $r$, control actions $\Phi$.
  \WHILE{system is (partially) observable and controllable}
  \STATE Estimate the system dynamics (I/O function) $\hat h$ using GP.
  \STATE Compute the best control action $u \in \Phi$ solving (\ref{e:contweigth1}).
  \STATE Compute variance, $v(y, u)$, of $\hat h$ as an estimate of $\mc I(\hat h)$.
  \STATE Update the data set $\mc D$ using newly observed data point $y$.
  \ENDWHILE
\end{algorithmic}
\end{algorithm}

\section{Examples}

\subsection{Dual Control of Logistic Map}

\begin{figure}[htp]
 \centering
 \includegraphics[width=0.8\columnwidth]{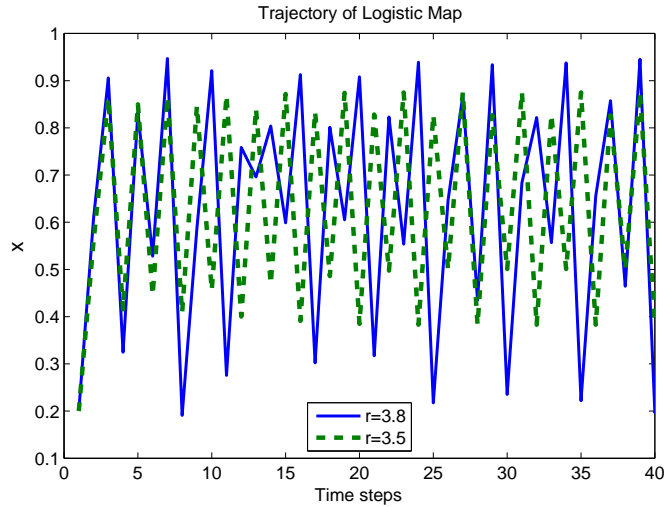}
 \caption{Example trajectories of the logistic map for $r=3.5$ and $r=3.8$.}
 \label{fig:logisticmap}
\end{figure}

The logistic map 
$$x(n+1)=r \, x(n) \left( 1- x(n)\right), $$
parameterized by the scalar $r$ is a well-known one-dimensional discrete-time nonlinear system, where $n$ denotes the time step or iteration. It is chosen as an illustrative example due to its interesting properties
%, which depend on the parameter $r$ as well as for visualization purposes. 
and for visualization purposes.
For $r=3.5$, logistic map converges to a limit cycle while it exhibits chaotic behavior for $r=3.8$ as shown in Figure~\ref{fig:logisticmap}.

\subsubsection*{Linear Control:}

First, the logistic map is controlled with additive actions while being identified using the GP method described in Algorithm~\ref{alg:algctrl1}:
$$ x(n+1)=r \, x(n) \left( 1- x(n)\right)+u(n).$$
The controller knows here that the control is linear (additive), and utilizes this extra knowledge in identifying
the system which simplifies the problem significantly. The system description (input-output relationship) from the perspective of the controller is:
$$ y(n+1)= \hat h(y(n))+u(n).$$
The control actions are taken from the finite set
$$\Phi=\{u_i \in [-1,1] : u_{i+1}=u_i+0.02, \; i=1,\ldots,101\}.$$ 

\begin{figure}[htp]
 \centering
 \includegraphics[width=0.8\columnwidth]{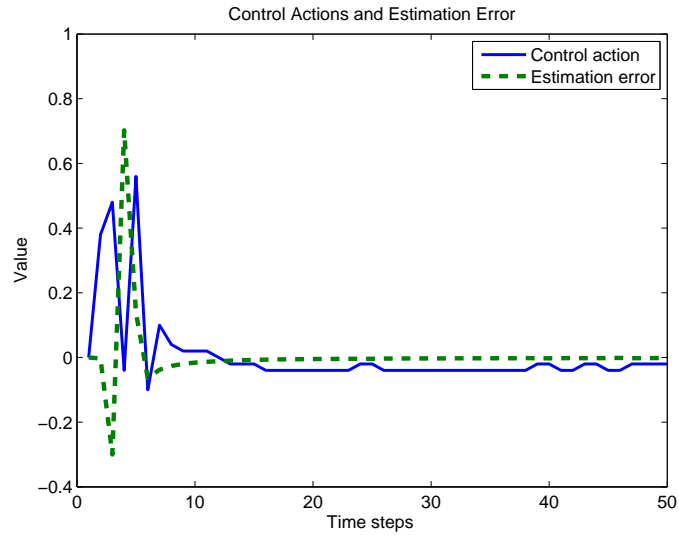}
 \caption{The control actions and state estimation errors for logistic map with $r=3.5$ and linear control.}
 \label{fig:control2}
\end{figure}

\begin{figure}[htp]
 \centering
 \includegraphics[width=0.8\columnwidth]{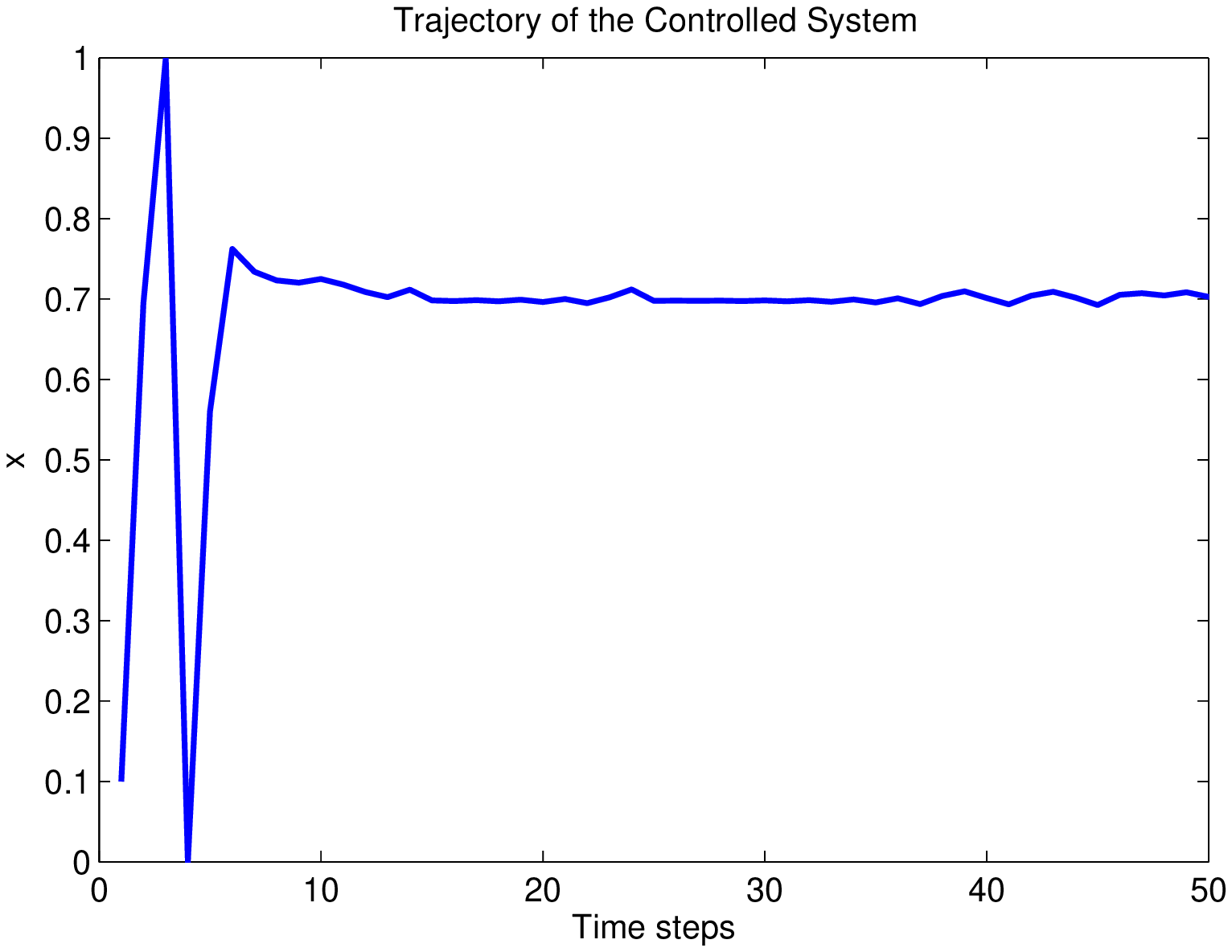}
 \caption{The controlled trajectory of the logistic map for $r=3.5$ and linear control.}
 \label{fig:trajectory2}
\end{figure}

\begin{figure}[htp]
 \centering
 \includegraphics[width=0.8\columnwidth]{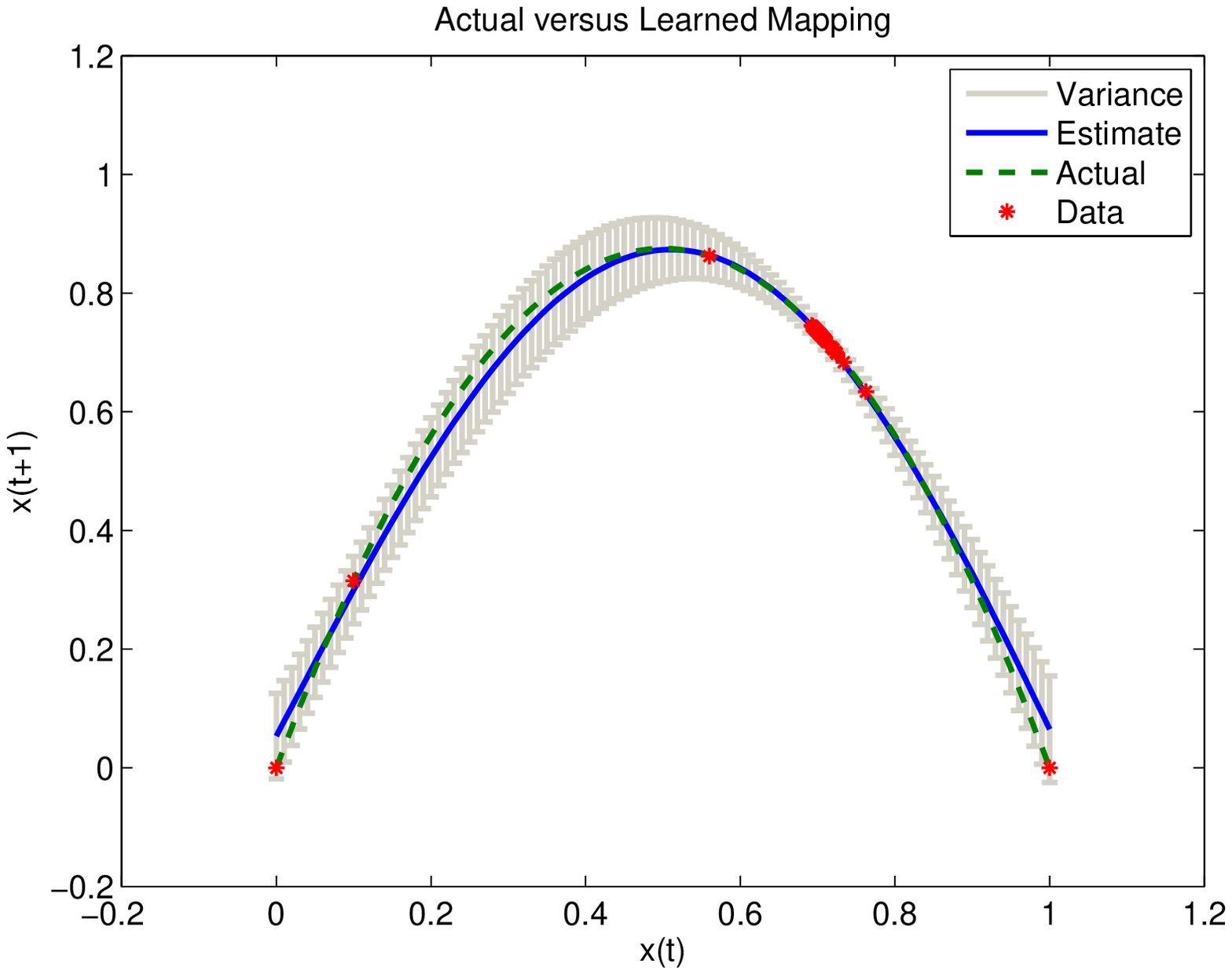}
 \caption{The logistic map and its estimate along with one standard deviation for $u=0$ and $r=3.5$ after $100$ iterations (data points).}
 \label{fig:mapping2}
\end{figure}

The kernel variance is $0.5$ and the weights in the objective function (\ref{e:contweigth1}) are chosen as $w_1=w_2=1$. The goal is stabilize the system at $x^*=0.8$, which constitutes the constant reference signal. The starting point is $x_0=0.1$. The control actions and state estimation errors over time (in each step based on arrived data points) for $r=3.5$ and the corresponding trajectory of the logistic map are depicted in Figures~\ref{fig:control2} and \ref{fig:trajectory2}. Note that, in this case the logistic map acts only as a nonlinear system with a limit cycle rather than behaving chaotically.
It is observed that approximately the first $10$ steps are used by the algorithm to explore or learn the system after which the trajectory approaches to the target. 
%The process involves learning up to as many input/output functions $h_u$ as the available control actions, i.e. the cardinality of $\Phi$.
The Figure~\ref{fig:mapping2} shows the estimated function versus the original mapping for $u=0$ as well as one standard deviation from estimated value. It can be seen that the variance is minimum, i.e. the estimate is best, around the target value.

%  However, the objective of the Algorithm~\ref{alg:algctrl1} is not learning the whole system but achieving the target in a greedy manner. Therefore, the system is estimated only partially. The Figure~\ref{fig:mapping2} shows the estimated function versus the original mapping for $u=0$ as well as one standard deviation from estimated value. It can be seen that the variance is minimum, i.e. the estimate is best, around the target value.


\begin{figure}[htp]
 \centering
 \includegraphics[width=0.8\columnwidth]{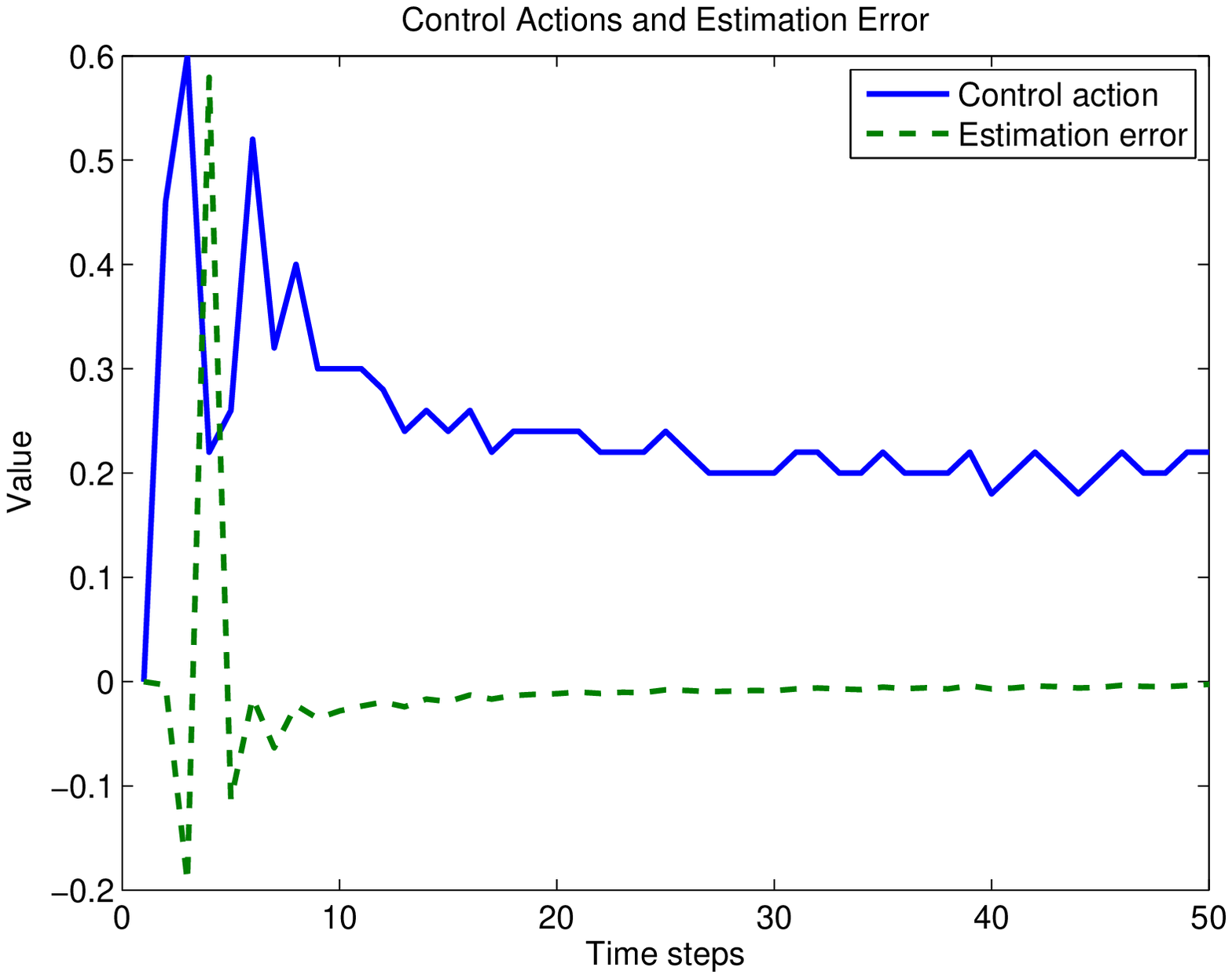}
 \caption{The control actions and state estimation errors for logistic map with $r=3.8$ and linear control}
 \label{fig:control1}
\end{figure}

\begin{figure}[htp]
 \centering
 \includegraphics[width=0.8\columnwidth]{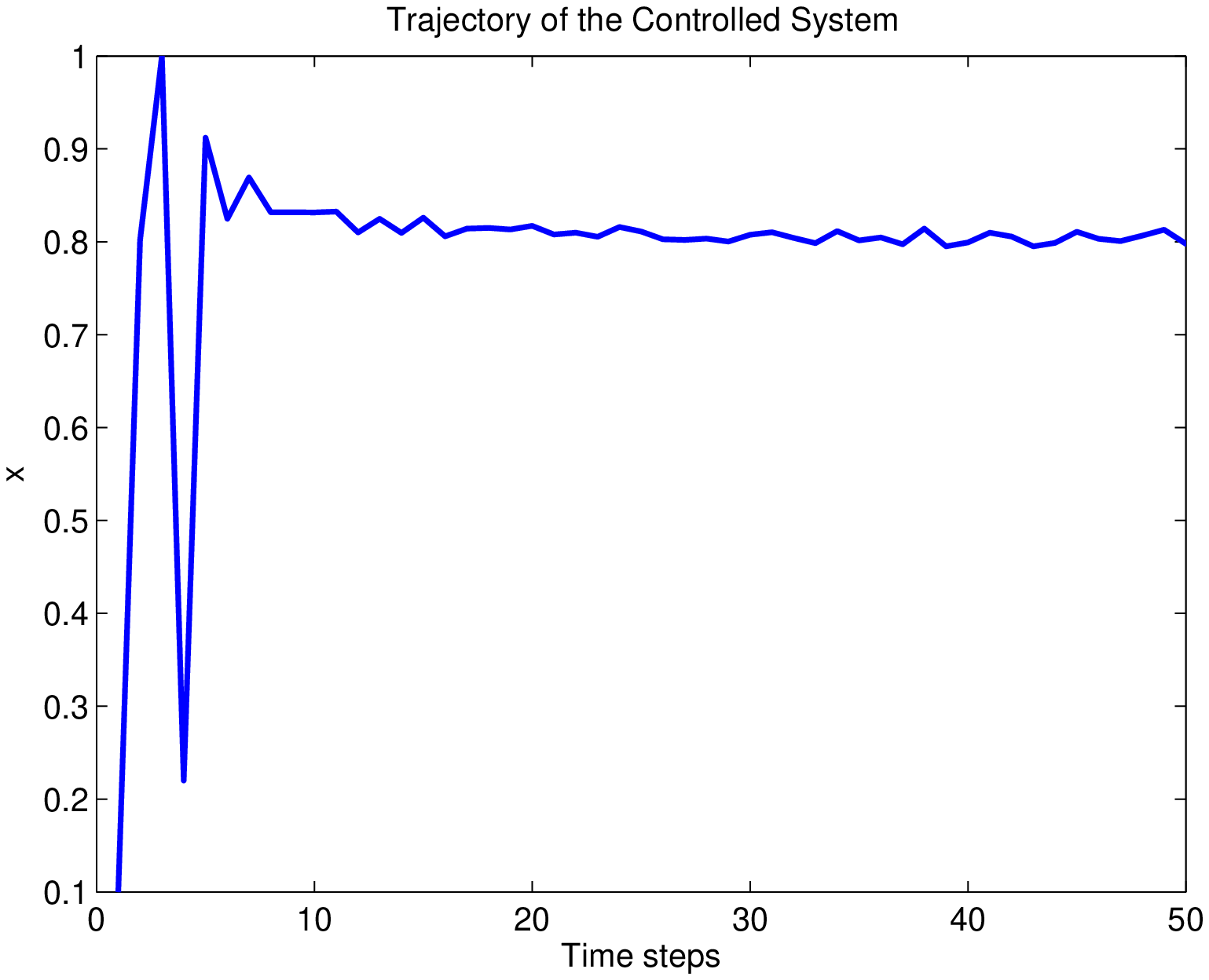}
 \caption{The controlled trajectory of the logistic map for $r=3.8$ and linear control.}
 \label{fig:trajectory1}
\end{figure}

\begin{figure}[htp]
 \centering
 \includegraphics[width=0.8\columnwidth]{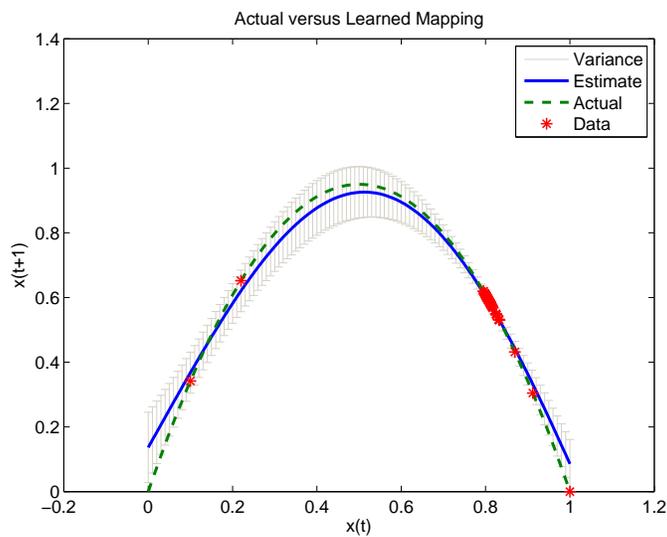}
 \caption{The logistic map and its estimate along with one standard deviation for $u=1.5$ and $r=3.8$ after $100$ iterations (data points).}
 \label{fig:mapping1}
\end{figure}

The same numerical analysis is repeated for $r=3.8$ in which case the logistic map behaves chaotically and the task turns to from control of an unknown nonlinear system to control of an unknown chaotic system. In this case, the goal is to stabilize the system at $x^*=0.8$. The control actions and state estimation errors over time (in each step based on arrived data points) for $r=3.8$ and the corresponding trajectory of the logistic map are depicted in Figures~\ref{fig:control1} and \ref{fig:trajectory1}. Note that the learning process takes longer in this case possibly due to the chaotic (complex) behavior of the system. The mapping shown in Figure~\ref{fig:mapping1}  shows the estimated function versus the original mapping for $u=1.5$. 

\subsubsection*{Nonlinear and Unknown Control:} Next, the logistic map is controlled with actions that affect the
system nonlinearly in a way that is unknown to the controller:
$$ x(n+1)=3.8 \, x(n) \left( 1- x(n)\right)+\cos(u).$$
The system description (input-output relationship) from the perspective of the controller is:
$$ y(n+1)= \hat h(y(n), u(n)).$$
Compared to the linear and known control case, this problem is obviously much harder to address. 
The control actions are taken from the finite set
$$\Phi=\{u_i \in [0,\pi] : u_{i+1}=u_i+0.1, \; i=1,\ldots,32\}.$$ 
The weights in the objective function (\ref{e:contweigth1}) are chosen initially as $w_1=1$ and $w_2=0$ to emphasize
exploration in the beginning but  $w_2$ is increased gradually to $w_2=40$ to achieve as good control performance as
possible.

\begin{figure}[htp]
 \centering
 \includegraphics[width=0.8\columnwidth]{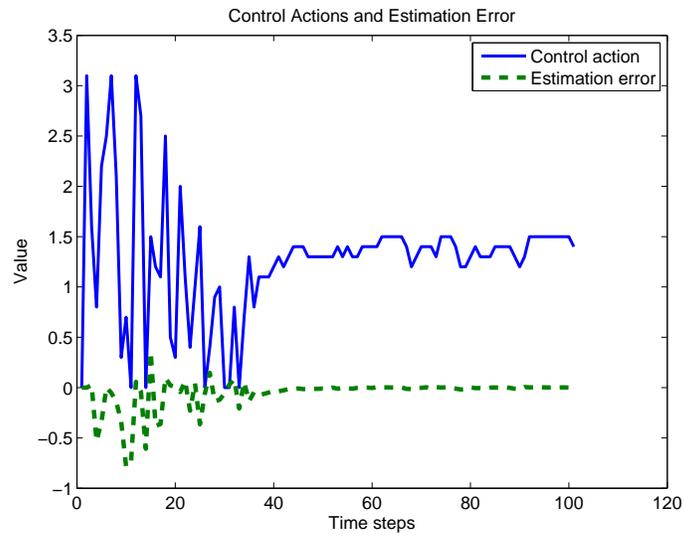}
 \caption{The control actions and state estimation errors for logistic map with $r=3.8$ and nonlinear control.}
 \label{fig:control3}
\end{figure}

\begin{figure}[htp]
 \centering
 \includegraphics[width=0.8\columnwidth]{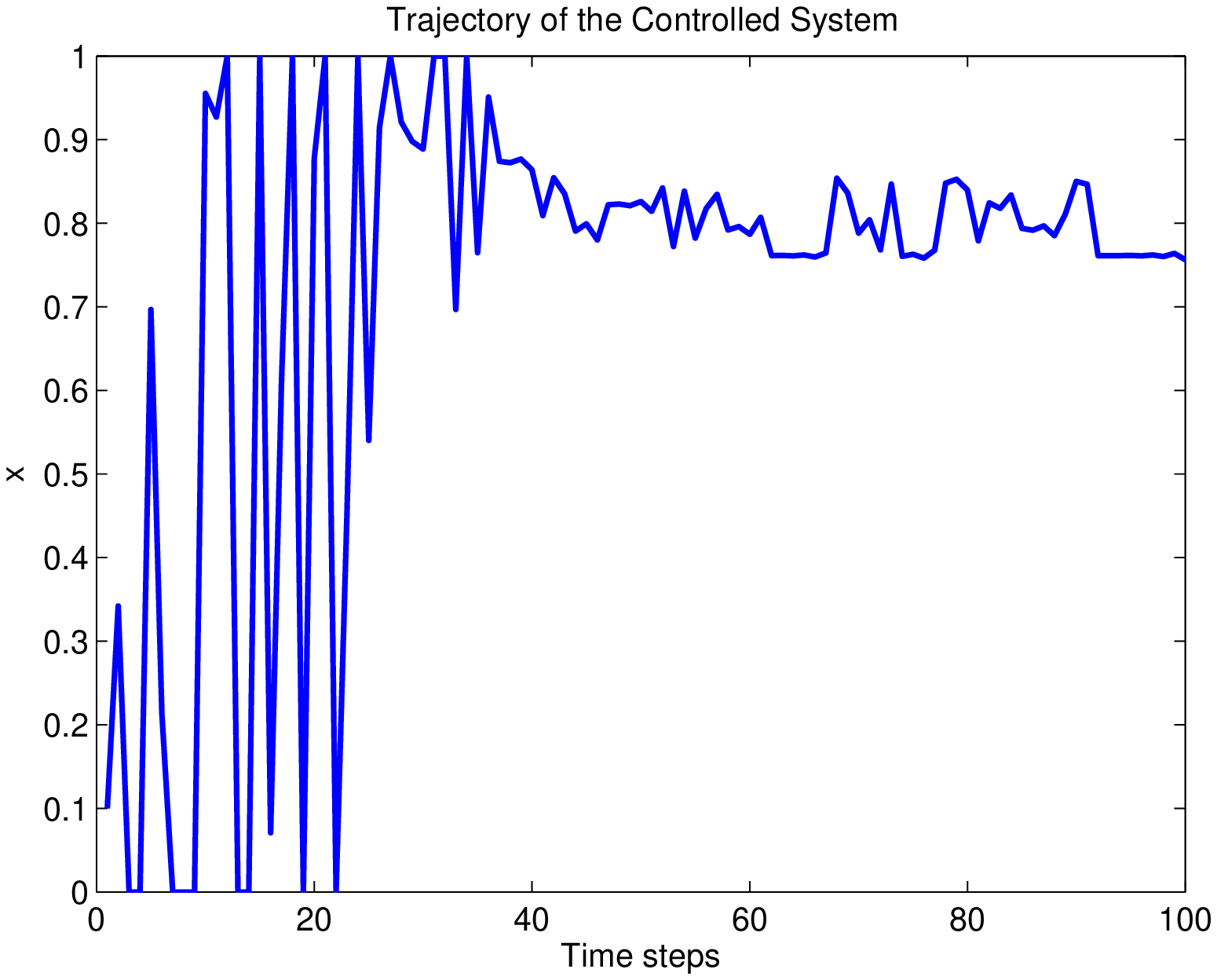}
 \caption{The controlled trajectory of the logistic map for $r=3.8$ and nonlinear control.}
 \label{fig:trajectory3}
\end{figure}

\begin{figure}[htp]
 \centering
 \includegraphics[width=0.8\columnwidth]{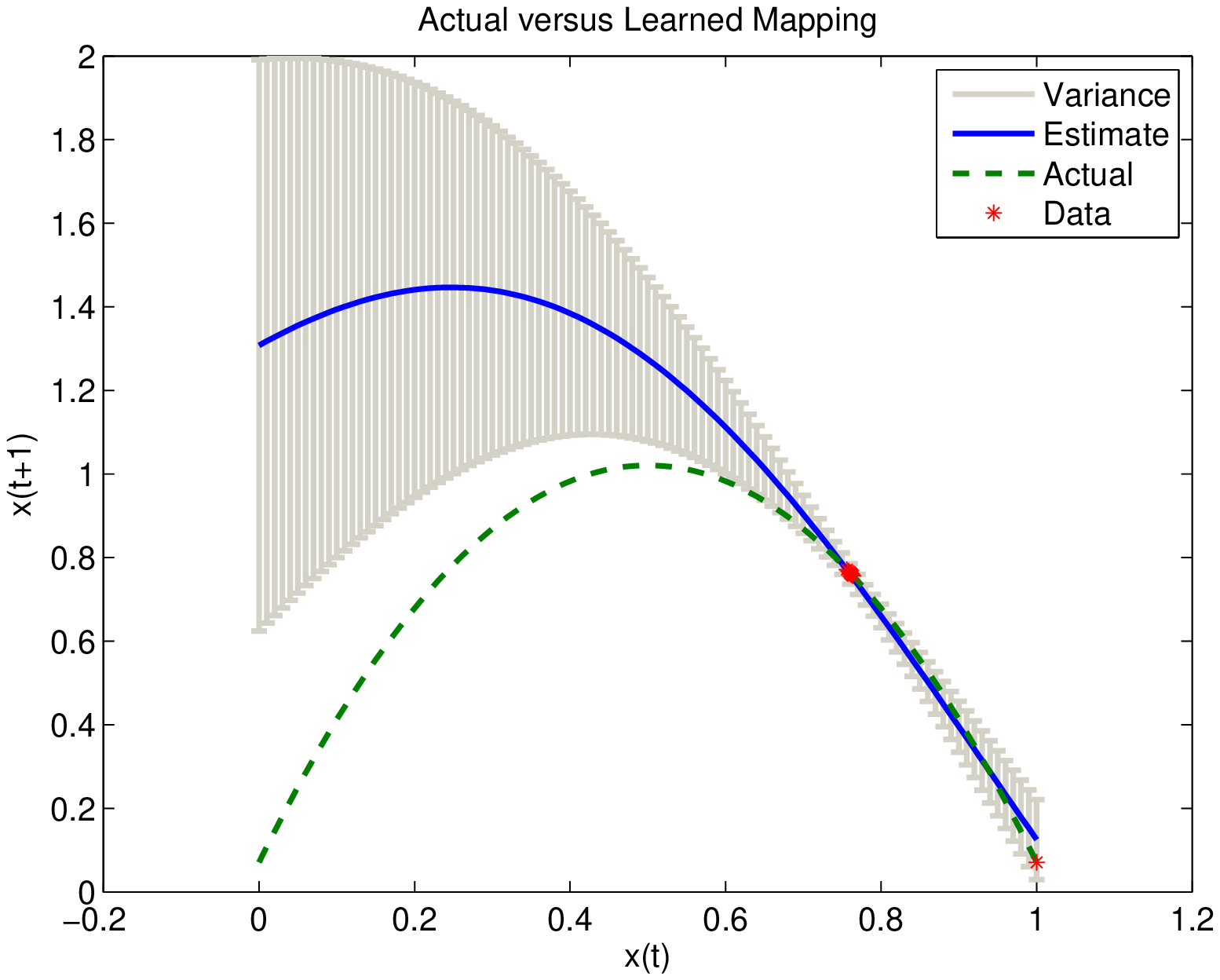}
 \caption{The logistic map and its estimate along with one standard deviation for $u=1.5$ and $r=3.8$ after $100$ iterations (data points).}
 \label{fig:mapping3}
\end{figure}

Figures~\ref{fig:control3}, \ref{fig:trajectory3}, and~\ref{fig:mapping3} summarize the obtained results. Since the objective of the Algorithm~\ref{alg:algctrl1} is not only learning the entire system behavior but achieving the control target in a greedy manner, the system is estimated accurately only around the target value. It is observed that the
learning process takes longer (twice as much of the case in the linear control) and the control actions are less accurate.
It should be kept in mind, however, that concurrently identifying and adaptively controlling a chaotic system with limited information is not an easy task.
 
\subsection{Position Control of a Cart with Inverted Pendulum}

The inverted pendulum on a cart is a classic example system for control problems. In this case, the problem is
formulated as the position control of the cart with the inverted pendulum, which is defined by the following
set of discrete-time nonlinear state-space equations \cite{wang1,wang2}:

\begin{align} 
 x_1(n+1)=x_1(n)+ T\, x_2(n), \quad \label{e:cart1}   
\end{align}
\begin{align}\label{e:cart2} 
 x_2(n+1)=x_2(n)+ \dfrac{T}{M+m\sin^2(x_3(n))}\left[ \mathbf{u}(n)  \right. \\ %\nonumber \\
 \nonumber  + m L x_4^2(n) \sin(x_3(n))-b  x_2(n)  \\
  \nonumber \left.- m g \cos(x_3(n))\sin(x_3(n))  \right] , 
\end{align}
\begin{align}
  x_3(n+1)=x_3(n)+ T\, x_4(n), \quad \label{e:cart3}\\
  x_4(n+1)=x_4(n)+\dfrac{T}{L\left( M+m \sin^2(x_3(n))\right) }  \label{e:cart4} \\
 \left[ - \mathbf{u}(n)\cos(x_3(n)) + (M+m) g\sin(x_3(n)) \right. \nonumber \\
 \nonumber \left. b x_2(n) \cos(x_3(n)) - m L x_4^2(n) \cos(x_3(n)) \sin(x_3(n)) \right],  \\
 y(n)=x_1(n), \quad \label{e:cart5}
\end{align}
where $T=0.05$ is the sampling period, $y=x_1$ is the position of the cart, $x_2=d x /dt$ is the cart velocity
$x_3=\theta$ is the inverted pendulum angle, $x_4=d \theta / d t$ is the angular velocity. The parameter
values are: $b=12.98$, $M=1.378$, $L=0.325$, $g=9.8$, and $m=0.051$. Further details on this standard model are
available in \cite{wang1,wang2}.

\begin{figure}[htp]
 \centering
 \includegraphics[width=0.8 \columnwidth]{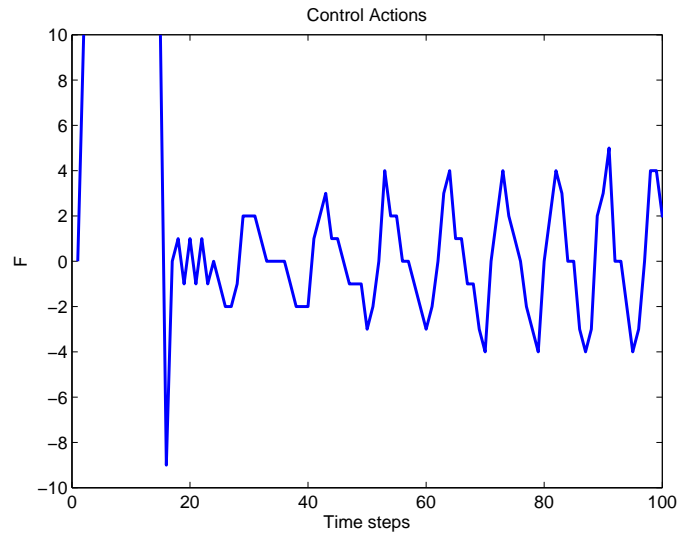}
 \caption{The control actions for the cart with inverted pendulum, single-step look ahead, and full knowledge.}
 \label{fig:invcontrol1}
\end{figure}

\begin{figure}[htp]
 \centering
 \includegraphics[width=0.8 \columnwidth]{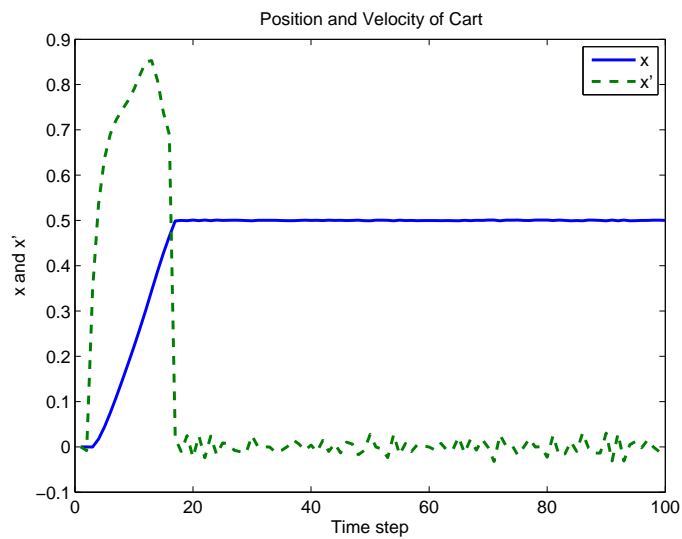}
 \caption{The trajectory of the cart with inverted pendulum controlled with full knowledge and single-step look ahead.}
 \label{fig:invtrajectory1}
\end{figure}

First, the cart is controlled using a one-step look-ahead strategy \textit{with full knowledge} from the starting point $x=[0,\, 0,\, 0.6,\, 0]$ with control actions chosen from the set $\{u_i \in [-10, 10] : u_{i+1}=u_i+1, \; i=1,\ldots,21\}$. The objective is to fix the position of the cart to $y^*=0.5$. The weights in the objective function (\ref{e:contweigth1}) are $w_1=1$ and $w_2=20$.
The results of this case shown in Figures~\ref{fig:invcontrol1} and \ref{fig:invtrajectory1} provide a benchmark to compare against.

\begin{figure}[htp]
 \centering
 \includegraphics[width=0.8 \columnwidth]{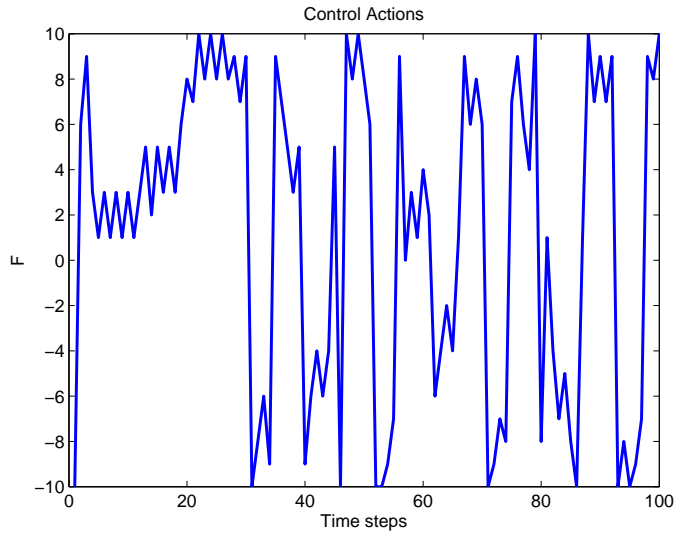}
 \caption{Dual control of the cart with inverted pendulum and single-step look ahead.}
 \label{fig:invcontrol2}
\end{figure}

\begin{figure}[htp]
 \centering
 \includegraphics[width=0.8 \columnwidth]{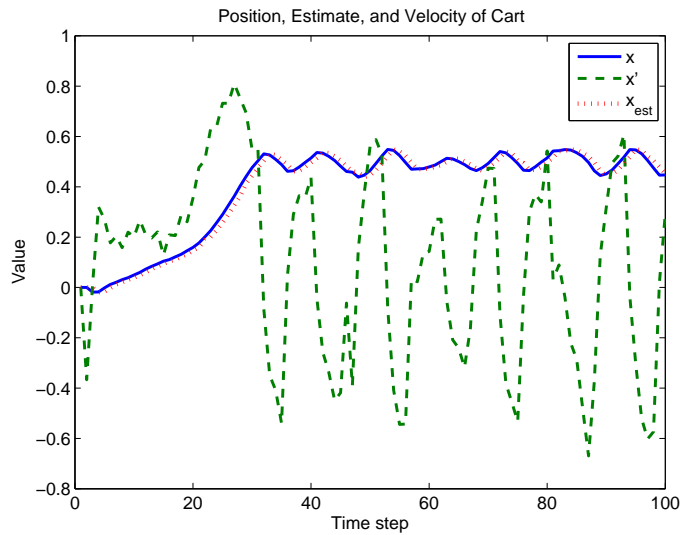}
 \caption{The trajectory of the cart with inverted pendulum under dual control with single-step look ahead estimates.}
 \label{fig:invtrajectory2}
\end{figure}

Next, the cart is controlled using a one-step look-ahead strategy as a \textit{as black-box system};
$y(n+1)=\hat h(\dot y(n), u(n))$. As side information, the
controller knows (\ref{e:cart1}), but has to estimate (\ref{e:cart2}) while (\ref{e:cart3}) and (\ref{e:cart4}) 
effectively act as external/unmodeled dynamics. The kernel and noise variance in GP are chosen as $0.5$ and $0.01$, respectively. The results obtained using Algorithm~\ref{alg:algctrl1} are shown in Figures~\ref{fig:invcontrol2} and \ref{fig:invtrajectory2}. The performance is satisfactory considering that the trajectory is within $10\%$ distance of the target within $30$ steps.


% -----------------------------------------------------------
\section{Literature Review} \label{sec:literature}

% Decision making with limited information is related to search theory. The idea of using information (theory) in this context is hardly new as evidenced by the article ``A New Look at the Relation Between Information Theory and Search Theory'' from 1979 \cite{pierce}. The subject is further studied in \cite{jaynes}. The topic of optimal search is more recently revisited by \cite{Zhu-search}, which contains substantial historical notes and studies problems where the search target distribution in itself is unobservable.

The book \cite{MacKaybook} provides important and valuable insights into the relationship between information theory, inference, and learning, where measuring information content of data points using Shannon information is discussed. However, focusing mainly on more traditional coding, communication, and machine learning topics, the book does not discuss the type of control problems presented in this paper. 

Learning plays an important role in the presented framework, especially \textit{regression}, which is a classical machine (or statistical) learning method. A very good introduction to the subject can be found in \cite{Bishopbook}. A complementary and detailed discussion on kernel methods is in \cite{schoelkopfbook}. Another relevant topic is Bayesian inference \cite{Tipping,MacKaybook}, which is in the foundation of the presented framework. In machine learning literature, Gaussian processes (GPs) are getting increasingly popular due to their various favorable characteristics. The book \cite{GPbook} presents a comprehensive treatment of GPs. Additional relevant works on the subject include \cite{MacKaybook,schoelkopfbook,MacKayGP}, which also discuss GP regression.

% Convex optimization \cite{boydbook} is a well-understood topic that is often easy to handle even if available information is limited. Optimizing nonconvex functions, however, is still a research subject \cite{globaloptsurvey}. It is interesting to note that the method known as \textit{kriging} in global optimization is almost the same as GP regression in machine learning. The field \textit{stochastic programming} focuses on optimization under uncertainty but assumes a certain amount of prior knowledge on the problem at hand and models the uncertainty probabilistically \cite{stochasticsurvery}. The popular heuristic method \textit{simulated annealing} \cite{simannealing} is essentially based on iterative random search. Another popular heuristic scheme particle swarm optimization \cite{swarmorig} is also based on random search but parallel in nature as a distinguishing characteristic rather than iterative. 

Gaussian processes have been recently applied to the area of optimization and regression \cite{BoylePhD} as well as system identification \cite{ThompsonPhD}. While the latter mentions active learning \cite{activelearning}, neither work discusses explicit information quantification or builds a connection with Shannon information theory. Using GP for system identification is discussed again in \cite{kocijan}, yet again without information collection aspects.
%The recent articles \cite{GPopt,turner1}, which utilize GP regression for optimization  in a setting similar to the one in this paper and for state-space inference and learning, respectively, do not consider information-theoretic aspects of the problem, either. Likewise,
% the article \cite{kriging1} on stochastic black box optimization, which considers a problem similar to the one here, does not take into account explicit measurement of information.
%The area of active learning or experiment design focuses on data scarcity in machine learning and makes use of Shannon information theory among other criteria \cite{activelearning}. 
The paper \cite{MacKaydataselect} discusses in a static optimization setting objective functions which measure the expected informativeness of candidate measurements within a Bayesian learning framework. 
The subsequent study \cite{gpactive1} investigates active learning for GP regression in machine learning applications using variance as a (heuristic) confidence measure  for test point rejection.
%Although this paper presents an approach similar to the one in this paper and interesting results, some of its comments and criticism are preliminary in nature and argued in an ad-hoc manner. 
%It is interesting that one of the final conclusions of  \cite{MacKaydataselect}, non-applicability of entropy as a measure of information in the search problem, is similar to those earlier works mentioned in the very first two articles cited in this section.

Dual control is an old topic, which has attracted the interest of the research community in the second half of the last century \cite{wittenmark1}. The article \cite{yame} revisits this subject and incorporates information explicitly into the dual control problem, but focuses on estimation of parameters in a known, linear system. Adopting a different perspective, a dynamic programming approach is presented recently in \cite{gpdynamic}, where 
an approximate value-function based reinforcement learning algorithm based on GPs and its online variant are presented. An application of GP-based identification and control to an autonomous blimp is discussed in \cite{gpblimp}.

% -----------------------------------------------------------
\section{Conclusion} \label{sec:conclusion}

The dual control approach presented in this paper addresses focuses on black-box control with very limited information.
The information acquired at each control step is quantified using the entropy measure from information theory and serves as the training input to a state-of-the-art Gaussian process regression (Bayesian learning) method. The quantification of the information obtained from each data point allows for iterative and joint optimization of both identification and control objectives. The results obtained from two illustrative examples, control of logistic map as a chaotic system and position control of a cart with inverted pendulum, demonstrate the developed approach.

The dynamic control problem in this paper differs from the static optimization analysis in \cite{valuetools11} in multiple ways. One of the main differences is the fact that the system states are now influenced indirectly through control actions. The data points used for identifying the underlying system mapping can only be selected indirectly (unlike static optimization) and under the constraints imposed by the nature of the ``control'' in the dynamic system at hand. 

The presented results should be considered mainly as an initial step. Future research directions are abundant
and include further investigation of the exploration-exploitation trade-off, more elaborate adaptive weighting parameters, and random sampling methods for problems in higher dimensional spaces. Applications to multi-person decision-making and game theory constitute another interesting future research topic.

% -----------------------------------------------------------
\section*{Acknowledgement}
This work is supported by Deutsche Telekom Laboratories.

%-------------------------------------------------------------------------
%\bibliographystyle{IEEEtranS}
%\bibliography{liminfocontrol}

% Generated by IEEEtranS.bst, version: 1.12 (2007/01/11)

\end{document}